\documentclass[11pt, reqno]{amsart}
\usepackage{amssymb,latexsym,amsmath,amsfonts,mathdots,enumitem}
\usepackage{latexsym}
\usepackage[mathscr]{eucal}
\usepackage{colortbl,xcolor}
\usepackage{lmodern}
\usepackage{sansmathaccent}
\usepackage[latin1]{inputenc}
\usepackage{tikz}
\usepackage{tikz-cd}
\usepackage{physics}
\usetikzlibrary{shapes,arrows}
\allowdisplaybreaks
\usetikzlibrary{matrix,calc,shapes,arrows,positioning}
\pdfmapfile{+sansmathaccent.map}

\numberwithin{equation}{section}
\theoremstyle{plain}

\newtheorem{thm}{Theorem}[section]

\newtheorem{cor}[thm]{Corollary}
\newtheorem{lem}[thm]{Lemma}
\newtheorem{prop}[thm]{Proposition}
\theoremstyle{definition}
\newtheorem{defn}[thm]{Definition}
\newtheorem{rem}[thm]{Remark}

\numberwithin{equation}{section}

\def\beq{\begin{eqnarray}}
	\def\eeq{\end{eqnarray}}
\def\beqa{\begin{eqnarray*}}
	\def\eeqa{\end{eqnarray*}}

\def\Ran{\operatorname{Ran}}

\def\beqn{\begin{equation}}
	\def\eeqn{\end{equation}}

\def\mg#1{}

\def\Ran{\operatorname{Ran}}

\renewcommand{\epsilon}{\varepsilon}
\renewcommand{\phi}{\varphi}

\renewcommand{\bf}[1]{\textbf{#1}}

\renewcommand{\sf}[1]{\textsf{#1}}

\DeclareMathOperator{\id}{\sf{id}}

\numberwithin{equation}{section}
\allowdisplaybreaks[4] 

\setlist[enumerate]{font=\upshape,noitemsep, topsep=0pt} 
\setlist[itemize]{noitemsep, topsep=0pt}

\begin{document}
	
	\title{A Note on Fundamental Operators Associated with $\mathbf{\Theta}_n$-Contraction}
	\author{Bhaskar Paul}
	
	\address[B. Paul]{Department of Mathematics, IIT Bhilai, 6th Lane Road, Jevra, Chhattisgarh 491002}
	\email{B. Paul:bhaskarpaul@iitbhilai.ac.in }
	
	\subjclass[2010]{47A13, 47A15, 47A20, 47A25, 47A45.}
	
	\keywords{Symmetrized polydisc, $\mathbf{\Theta}_n$-contraction, fundamental operators, unitary invariant, complete spectral set}
	
	\maketitle

	\begin{abstract}
		An $n$-tuple of commuting bounded operators
		$\mathbf{T} = (T_1, \dots, T_n)$ on a Hilbert space $\mathcal{H}$ is called a $\mathbf{\Theta}_n$-contraction
		if $\overline{\mathbf{\Theta}}_n$ is a spectral set for $\mathbf{T}$. The fundamental equations of $\mathbf{T}$
		is defined by
		\begin{equation*}
			\begin{aligned}
				T_i-T_{n-i}^*T_n^p
				&= \sum_{l=0}^{p-1}
				T_n^{*p-1-l}D_{T_n}A^{(i)}_l
				D_{T_n}T_n^{p-1} + \sum_{l=0}^{p-2} T_n^{*p-2-l}D_{T_n}A^{(i)}_l
				D_{T_n}T_n^{p-2} + \dots\\
				&+ \sum_{l=0}^{1}
				T_n^{*1-l}D_{T_n}A^{(i)}_l
				D_{T_n}T_n + D_{T_n}A^{(i)}_0D_{T_n},
				\quad 1 \le i \le n-1,
			\end{aligned}
		\end{equation*}
		where $(A^{(i)}_0, \dots, A^{(i)}_p),$
		$1 \le i \le n-1$ are referred to as
		fundamental operators of $\mathbf{T}$.
		By using the non-uniqueness of tuples of fundamental operators, we show that there are $\mathbf{\Theta}_n$-contractions that possess $\mathbf{\Theta}_n$-isometric dilation, but the fundamental operators do not satisfy some of the conditions assumed in conditional dilation of
		$\mathbf{\Theta}_n$-contraction. We also prove that the tuples of fundamental operators do not satisfy
		\begin{equation*}
			\begin{aligned}
				D_{T_n}T_i &=
				\sum_{l=0}^{p} A^{(i)}_lD_{T_n}T^l_n,
				\quad 1 \le i \le n-1,
			\end{aligned}
		\end{equation*}
		which negatively answers the open question proposed in
		the paper
		[\textit{https://arxiv.org/abs/2608.03574v2}]. Finally, we demonstrate under certain conditions that for a given set of operators $\{A^{(i)}_l : 1 \le i \le n-1, 0 \le l \le p-1\}$ defined on a Hilbert space $\mathcal{E}$, there exists a $\mathbf{\Theta}_n$-contraction $\mathbf{T}$ such that $(A^{(i)}_0, \dots, A^{(i)}_{p-1})$,
		$1 \le i \le n-1$ are the fundamental operators of $\mathbf{T}$.
	\end{abstract}

	\section{Introduction}\label{Intro}
	
	We denote the open and closed unit disc in the complex plane by $\mathbb{D}$ and $\overline{\mathbb{D}}$ respectively. The \textit{symmetrization map} $\mathbf{s} = (s_1, \dots, s_n) : \mathbb{C}^n \to \mathbb{C}^n$
	is defined by
	\begin{equation*}
		\begin{aligned}
			s_i(z_1, \dots, z_n) &=
			\sum_{1 \le k_1 < \cdots < k_i \le n}
			z_{k_1} \cdots z_{k_i}, \quad 1 \le i \le n.
		\end{aligned}
	\end{equation*}
	Its image $\mathbb{G}_n = \mathbf{s}(\mathbb{D}^n)$ is the symmetrized polydisc, while $\Gamma_n = \mathbf{s}(\overline{\mathbb{D}}^n)$
	denotes its closure. A detailed study on $\mathbb{G}_n$ (respectively, $\Gamma_n$) and $\Gamma_n$-contraction can be found in \cite{Biswas, Costara}. S. Biswas, G. Ghosh, E. K. Narayanan, and S. Shyam Roy generalized the symmetrization map in \cite{Ghosh}. Let $p$ be a fixed divisor of $m$, let $\theta_0 = 1$, and define
	\begin{equation*}
		\begin{aligned}
			\theta_i(z_1, \dots, z_n) &=
			s_i(z^m_1, \dots, z^m_n), \quad
			1 \le i \le n-1,
		\end{aligned}
	\end{equation*}
	and
	\begin{equation*}
		\begin{aligned}
			\theta_n(z_1, \dots, z_n) &=
			(z_1 \cdots z_n)^q, \quad q = \frac{m}{p}.
		\end{aligned}
	\end{equation*}
	The images of $\mathbb D^n$ and $\overline{\mathbb{D}}^n$ under $\theta = (\theta_1, \dots, \theta_n)$ are denoted by $\mathbf{\Theta}_n$ and $\overline{\mathbf{\Theta}}_n$, respectively. Observe that $\mathbf{\Theta}_1 = \mathbb{D}$,
	and that $\mathbf{\Theta}_n = \mathbb{G}_n$ when $m = p = 1$. From \cite{Biswas 2, Ghosh} it is immediate that $(\theta_1, \dots, \theta_n) \in \overline{\mathbf{\Theta}}_n$ if and only if every zero of
	\begin{equation}\label{P}
		\begin{aligned}
			P(z) &= z^n - \theta_1z^{n-1} + \dots + (-1)^n\theta_n^p
		\end{aligned}
	\end{equation}
	lies in $\overline{\mathbb{D}}$. The characterization of domain $\mathbf{\Theta}_n$, and the interpolation by rational $\mathbf{\Theta}_n$-inner function is done in \cite{Paul}. Throughout this paper, we consider $m > 1, p > 1$ unless otherwise specified.
	
	Let $\Omega$ be a compact subset in $\mathbb{C}^n$, and let $\mathcal{O}(\Omega)$ denote the algebra of bounded holomorphic functions in some neighbourhood of $\Omega$, equipped with the sup norm $||\, . \,||_{\infty, \Omega}$, defined by
	\begin{equation*}
		\begin{aligned}
			||f||_{\infty, \Omega} &= \sup\{|f(z)| : z \in \Omega\}.
		\end{aligned}
	\end{equation*}
	A compact subset $\Omega$ of $\mathbb{C}^n$ is said to be a \textit{spectral set} for a commuting $n$-tuple of bounded operators $\mathbf{T} = (T_1, \dots, T_n)$ if $\sigma(\mathbf{T}) \subset \Omega$ and the homomorphism $\rho_{\mathbf{T}} : \mathcal{O}(\Omega) \to \mathcal{B}(\mathcal{H})$ defined by
	\begin{equation*}
		\begin{aligned}
			\rho_{\mathbf{T}}(1) &= I_{\mathcal{H}} \quad \text{and} \quad \rho_{\mathbf{T}}(z_i) = T_i, \quad 1 \le i \le n
		\end{aligned}
	\end{equation*}
	is contractive (see \cite{McCarthy 1, Nagy}). The notion of a spectral set for an operator on a Hilbert space was introduced by von Neumann. The following theorem says that the closed unit disc is a spectral set for every contraction.
	\begin{thm}[{\cite[Chapter 1, Corollary 1.2]{Nagy}}]\label{von Neumann Inequality}
		Let $T$ be a contraction on a
		Hilbert space. Then for any polynomial $p$,
		\begin{equation*}
			\begin{aligned}
				||p(T)|| &\le \sup\{|p(z)| : |z| \le 1\}.
			\end{aligned}
		\end{equation*}
	\end{thm}
	The following dilation theorem for a contraction is due to Sz.-Nagy.
	\begin{thm}[{\cite[Theorem 1.1]{paulsen}}]\label{Nagy Dilation}
		Let $T$ be a contraction operator on Hilbert space $\mathcal{H}$.  Then there is a Hilbert space $\mathcal{K}$ containing $\mathcal{H}$ as a subspace and a unitary operator $U$ on $\mathcal{K}$ such that
		\begin{equation}\label{Power Dilation}
			\begin{aligned}
				T^n &= P_{\mathcal{H}}U^n|_{\mathcal{H}}, \quad n \in \mathbb{N} \cup \{0\}.
			\end{aligned}
		\end{equation}
	\end{thm}
	The dilation presented in Theorem \ref{Nagy Dilation} is called \textit{power dilation} of $T$. Sch\"{a}ffer explicitly constructed an unitary dilation of a given contraction (see \cite{Schaffer}).
	
	Suppose $F = ((f_{ij}))$ is a matrix-valued polynomial defined on $\Omega$. Then $\Omega$ is said to be a \textit{complete spectral set} for $\mathbf T$ if
	\begin{equation*}
		\begin{aligned}
			||F(\mathbf{T})|| &\le ||F||_{\infty, \Omega}
		\end{aligned}
	\end{equation*}
	for every $F \in \mathcal{O}(\Omega) \otimes M_{k \times k}(\mathbb{C}), k \ge 1$. A commuting $n$-tuple $\mathbf{T}$ for which $\Omega$ is a spectral set is said to admit a $\partial \Omega$-normal dilation if there exist a Hilbert space $\mathcal{K} \supseteq \mathcal{H}$ and a commuting $n$-tuple of normal operators $\mathbf{N} = (N_1, \dots, N_n)$
	on $\mathcal{K}$ satisfying $\sigma(\mathbf{N}) \subseteq \partial \Omega$ such that
	\begin{equation*}
		\begin{aligned}
			P_{\mathcal H}f(\mathbf{N})|_{\mathcal H} &= f(\mathbf{T}), \quad f \in \mathcal{O}(\Omega).
		\end{aligned}
	\end{equation*}
	The following theorem due to Arveson (see \cite[Corollary, p. 279]{Arveson}) connects the notion of complete spectral set and $\partial \Omega$-normal dilation.
	\begin{thm}\label{Thm Arveson}
		Let $\mathbf{T} = (T_1, \dots, T_n)$ have K as a spectral set. Then $\mathbf{T}$ has $K$ as a complete spectral set if and only if $\mathbf{T}$ has a $\partial K$-normal dilation.
	\end{thm}
	
	Agler in \cite[Theorem 1.9]{Agler} proved the following result related to the complete contractivity of a unital algebra homomorphism from $\mathcal{O}(K)$ to $M_2(\mathbb{C})$, which we recall below. It plays a crucial role in proving one of our main results.
	
	\begin{thm}\label{Thm Agler}
		Let $K$ be a compact subset of $\mathbb{C}^n$. A unital $(\Pi(1) = 1)$ contractive $(||\Pi(f)|| \le ||f||)$ algebra homomorphism $\Pi : \mathcal{O}(K) \to M_2(\mathbb{C})$ is completely contractive.
	\end{thm}
	
	A functional model for a pure contraction was demonstrated by Sz.-Nagy and Foias in \cite{Nagy}. The development is briefed as follows. Let $T$ be a contraction a Hilbert space $\mathcal{H}$. We first recall that for a contraction $T$ on a Hilbert space $\mathcal{H}$, the defect operator and defect space are defined by
	\begin{equation*}
		\begin{aligned}
			D_T = (I - T^*T)^{1/2} \quad \text{and} \quad \mathcal{D}_T = \overline{\operatorname{Ran}}D_T,
		\end{aligned}
	\end{equation*}
	respectively. Then the operators $D_T$ and $D_{T^*}$ satisfy the following identity:
	\begin{equation*}
		\begin{aligned}
			TD_T &= D_{T^*}T \quad \text{equivalently} \quad D_TT^* = T^*D_{T^*}.
		\end{aligned}
	\end{equation*}
	The \textit{characteristic function} $\Theta_T$ of $T$ is defined as
	\begin{equation}\label{Characteristic}
		\begin{aligned}
			\Theta_T(z)
			&= (- T + D_{T^*}(I - zT^*)^{-1}D_T)|_{\mathcal{D}_T}, \quad z \in \mathbb{D}.
		\end{aligned}
	\end{equation}
	Note that $\Theta_T \in \mathcal{B}(\mathcal{D}_T, \mathcal{D}_{T^*})$.  We define the  multiplication operator $M_{\Theta_T} : H^2(\mathbb{D}) \otimes \mathcal{D}_T \to H^2(\mathbb{D}) \otimes \mathcal{D}_{T^*}$ by
	\begin{equation*}
		\begin{aligned}
			M_{\Theta_T}f(z) = \Theta_T(z)f(z), \quad z \in \mathbb{D}.
		\end{aligned}
	\end{equation*}
	Consider the quotient space $\mathcal{H}_T = (H^2(\mathbb{D}) \otimes \mathcal{D}_{T^*}) \ominus M_{\Theta_T}(H^2(\mathbb{D}) \otimes \mathcal{D}_T)$. The space $\mathcal{H}_T$ is called the \textit{model space for $T$}. We now recall the functional model for pure contraction from \cite{Nagy}.
	
	\begin{thm}\label{Pure Contraction Model}
		Every pure contraction $T$ defined on a Hilbert space $\mathcal{H}$ is unitarily 
		equivalent to the operator $T_1$ on the Hilbert space $\mathcal{H}_T = (H^2(\mathbb{D}) \otimes \mathcal{D}_{T^*}) \ominus M_{\Theta_T}(H^2(\mathbb{D}) \otimes \mathcal{D}_T)$ defined as
		\begin{equation*}
			\begin{aligned}
				T_1 &= \mathbb{P}_{\mathcal{H}_T}
				(M_z \otimes I_{\mathcal{D}_{T^*}})|_{\mathcal{H}_T}.
			\end{aligned}
		\end{equation*}
	\end{thm}
	
	During the last two decades, many authors have developed functional models for pure $\Gamma$-contraction, pure tetrablock contraction, and pure $\Gamma_n$-contraction in \cite{S. Pal 2, S. Pal 3, Sau}. A functional model for a class of pure $\mathbf{\Theta}_n$-contraction can also be found in \cite[Theorem 3.6]{Paul 2}.
	
	Let us recall some definitions associated with operator theory on $\mathbf{\Theta}_n$ from \cite{Biswas 2}.
	
	\begin{defn}\label{Defition 1}
		\begin{enumerate}
			\item Let $\mathbf{N} = (N_1, \dots, N_n)$ be a commuting $n$-tuple of bounded normal operators on a Hilbert space $\mathcal H$. We say that $\mathbf{N}$ is a \textit{$\mathbf{\Theta}_n$-unitary} if its joint spectrum is contained in the distinguished boundary $b\mathbf{\Theta}_n$.
			
			\item Let $\mathbf{V} = (V_1, \dots, V_n)$ be a commuting $n$-tuple of bounded operators on a Hilbert space $\mathcal H$. We say that $\mathbf{V}$ is a \emph{$\mathbf{\Theta}_n$-isometry} if there exist a Hilbert space $\mathcal{K} \supseteq \mathcal{H}$ and a $\mathbf{\Theta}_n$-unitary $\mathbf{N} = (N_1, \dots, N_n)$ on $\mathcal K$ such that
			$V_i = N_i|_{\mathcal{H}}, 1\le i\le n.$
			We say that $\mathbf{V}$ is a \emph{$\mathbf{\Theta}_n$-co-isometry} if
			$\mathbf{V}^* = (V^*_1, \dots, V^*_n)$
			is a $\mathbf{\Theta}_n$-isometry. A $\mathbf{\Theta}_n$-isometry $\mathbf{V}$ is said to be a \textit{pure $\mathbf{\Theta}_n$-isometry} if $V_n$ is a pure isometry.
		\end{enumerate}
	\end{defn}
	
	A straightforward application of Theorem \ref{Thm Agler} is the following theorem. Since the existence of a $\mathbf{\Theta}_n$-isometric dilation ensures the existence of a $\mathbf{\Theta}_n$-unitary dilation, we can state the following theorem for a $\mathbf{\Theta}_n$-isometric dilation.
	
	\begin{thm}\label{2 times 2 Theta_n-Contraction}
		Every $2 \times 2$ matrix $\mathbf{\Theta}_n$-contraction has a $\mathbf{\Theta}_n$-isometric dilation.
	\end{thm}
	
	A major portion of the study of $\Gamma_n$-contractions relies on the theory of \textit{fundamental equations} and \textit{fundamental operators}. For a $\Gamma_n$-contraction $(S_1, \dots, S_n)$, the fundamental equations are 
	\begin{equation}\label{Gamma_n Fundamental}  
		\begin{aligned}
			S_i - S_{n-i}^*S_n = D_{S_n}E_iD_{S_n}, \quad S_{n-i} - S_i^*S_n = D_{S_n}E_{n-i}D_{S_n}, \quad E_i, E_{n-i}\in\mathcal{B}(\mathcal{D}_{S_n}),
		\end{aligned}
	\end{equation} 
	for $1 \le i \le n-1$. The detailed study on this direction can be found in \cite{S. Pal, A. Pal}. It is natural to ask whether an analogous theory exists for $\mathbf{\Theta}_n$-contractions. Thus, all these developments motivate us to introduce the system of equations associated with a $\mathbf{\Theta}_n$-contraction $\mathbf{T} = (T_1, \dots, T_n)$ as follows:
	\begin{equation}\label{Fundamental}
		\begin{aligned}
			T_i - T^*_{n-i}T^p_n
			&= \sum_{l=0}^{p-1} T^{*p-1-l}_n D_{T_n}A^{(i)}_lD_{T_n} T^{p-1}_n + \sum_{l=0}^{p-2} T^{*p-2-l}_n D_{T_n}A^{(i)}_lD_{T_n} T^{p-2}_n + \dots\\
			& \quad + \sum_{l=0}^{1} T^{*1-l}_n D_{T_n}A^{(i)}_lD_{T_n} T_n + D_{T_n}A^{(i)}_0D_{T_n}, \quad 1 \le i \le n-1.
		\end{aligned}
	\end{equation}
	We refer to \eqref{Fundamental} as the \textit{fundamental equations} of $\mathbf{T}$, and
	\[(A^{(i)}_0, \dots, A^{(i)}_{p-1}), \quad 1 \le i \le n-1\]
	are the \textit{tuples of fundamental operators} of $\mathbf{T}$. Note that for $m = p = 1$, these equations reduce to \eqref{Gamma_n Fundamental} for $\Gamma_n$-contractions with $A^{(i)}_0 = E_i, 1 \le i \le n-1$. A. Pal in \cite[Theorem 4.4]{A. Pal} proved that the existence and uniqueness of the fundamental operators for $\Gamma_n$-contractions.
	
	We recall the conditional dilation for a given $\mathbf{\Theta}_n$-contraction in the following theorem.
	
	\begin{thm}[{\cite[Theorem 5.2]{Paul 1}}]\label{Conditional Dilation}
		Let $\mathbf{T} = (T_1, \dots, T_n)$ be a $\mathbf{\Theta}_n$-contraction on a Hilbert space $\mathcal H$, and let $A^{(i)}_l, 0 \le l \le p, 1 \le i \le n-1$, be bounded operators on $\mathcal D_{T_n}$. Let
		\begin{equation*}
			\begin{aligned}
				\mathcal{K} &= \mathcal{H} \oplus
				\mathcal{D}_{T_n} \oplus
				\mathcal{D}_{T_n}
				\oplus \cdots =
				\mathcal{H} \oplus \ell^2(\mathcal{D}_{T_n}).
			\end{aligned}
		\end{equation*}
		Define an $n$-tuple of operators
		$\mathbf{V} = (V_1, \dots, V_n)$
		on $\mathcal{K}$ by	
		\begin{equation}\label{V_i}
			\begin{aligned}
				V_i &=
				\left[
				\begin{smallmatrix}
					T_i &\hspace{0.5cm} 0 &\hspace{0.5cm} 0 &\hspace{0.5cm} 0 &\hspace{0.5cm} \dots \hspace{0.5cm}\\
					\displaystyle\sum_{l=0}^{p-1} A^{(i)}_{p-l}D_{T_n}T^{p-1-l}_n &\hspace{0.5cm} A^{(i)}_0 &\hspace{0.5cm} 0 &\hspace{0.5cm} 0 &\hspace{0.5cm} \dots \hspace{0.5cm}\\
					\displaystyle\sum_{l=0}^{p-2} A^{(i)}_{p-l}D_{T_n}T^{p-2-l}_n &\hspace{0.5cm} A^{(i)}_1 &\hspace{0.5cm} A^{(i)}_0 &\hspace{0.5cm} 0 &\hspace{0.5cm} \dots \hspace{0.5cm}\\
					\vdots & \hspace{0.5cm} \vdots & \hspace{0.5cm} \vdots & \hspace{0.5cm} \vdots & \hspace{0.5cm} \ddots \hspace{0.5cm}\\
					\displaystyle\sum_{l=0}^{1} A^{(i)}_{p-l}D_{T_n}T^{1-l}_n &\hspace{0.5cm} A^{(i)}_{p-2} &\hspace{0.5cm} A^{(i)}_{p-3} &\hspace{0.5cm} A^{(i)}_{p-4} &\hspace{0.5cm} \dots \hspace{0.5cm}\\
					A^{(i)}_pD_{T_n} &\hspace{0.5cm} A^{(i)}_{p-1} &\hspace{0.5cm} A^{(i)}_{p-2} &\hspace{0.5cm} A^{(i)}_{p-3} &\hspace{0.5cm} \dots \hspace{0.5cm}\\
					0 &\hspace{0.5cm} A^{(i)}_p &\hspace{0.5cm} A^{(i)}_{p-1} &\hspace{0.5cm} A^{(i)}_{p-2} &\hspace{0.5cm} \dots \hspace{0.5cm}\\
					\vdots &\hspace{0.5cm} \vdots &\hspace{0.5cm} \vdots &\hspace{0.5cm} \vdots &\hspace{0.5cm} \ddots \hspace{0.5cm}
				\end{smallmatrix}\right]
			\end{aligned}
		\end{equation}
		for $1 \le i \le n-1$ and
		\begin{equation}\label{V_n}
			\begin{aligned}
				V_n &=
				\left[
				\begin{smallmatrix}
					T_n &\hspace{0.5cm} 0 &\hspace{0.5cm} 0 &\hspace{0.5cm} 0 &\hspace{0.5cm} \dots \hspace{0.5cm}\\
					D_{T_n} &\hspace{0.5cm} 0 &\hspace{0.5cm} 0 &\hspace{0.5cm} 0 &\hspace{0.5cm} \dots \hspace{0.5cm}\\
					0 &\hspace{0.5cm} I &\hspace{0.5cm} 0 &\hspace{0.5cm} 0 &\hspace{0.5cm} \dots \hspace{0.5cm}\\
					0 &\hspace{0.5cm} 0 &\hspace{0.5cm} I &\hspace{0.5cm} 0 &\hspace{0.5cm} \dots \hspace{0.5cm}\\
					0 &\hspace{0.5cm} 0 &\hspace{0.5cm} 0 &\hspace{0.5cm} I &\hspace{0.5cm} \dots \hspace{0.5cm}\\
					\vdots &\hspace{0.5cm} \vdots &\hspace{0.5cm} \vdots &\hspace{0.5cm} \vdots &\hspace{0.5cm} \ddots \hspace{0.5cm}
				\end{smallmatrix}\right],
			\end{aligned}
		\end{equation}
		where $I$ denotes the identity operator on $\mathcal{D}_{T_n}$. Then $\mathbf{V} = (V_1, \dots, V_n)$ is a $\mathbf{\Theta}_n$-isometric dilation of $\mathbf{T}$ provided that $(\gamma_1V_1, \dots, \gamma_{n-1}V_{n-1})$ is a $\Gamma_{n-1}$-contraction and the following conditions are satisfied:
		\begin{enumerate}
			\item[(1)] $A^{(i)*}_l = A^{(n-i)}_{p-l}, 1 \le i \le n-1, 0 \le l \le p$, and $A^{(i)}_l = 0$ for $l>p$.
			
			\item[(2)] $\displaystyle\sum_{l=0}^{k} [A^{(i)}_l, A^{(j)}_{k-l}]=0$, $0 \le l \le 2p$ and $1 \le i, j \le n-1$, where $[A, B] = AB - BA$.
			
			\item[(3)] $D_{T_n}T_i = \displaystyle\sum_{k=0}^{p} A^{(i)}_kD_{T_n}T_n^k, 1 \le i \le n-1$.
			
			\item[(4)] $T_i-T_{n-i}^*T_n^p
			= \displaystyle\sum_{l=0}^{p-1}T^{*p-1-l}_nD_{T_n}A^{(i)}_lD_{T_n}T^{p-1}_n + \displaystyle\sum_{l=0}^{p-2}T_n^{*p-2-l}D_{T_n}A^{(i)}_lD_{T_n}T^{p-2}_n
			+ \dots\\
			+ \displaystyle\sum_{l=0}^{1}T^{*1-l}_nD_{T_n}A^{(i)}_lD_{T_n}T_n + D_{T_n}A^{(i)}_0D_{T_n}$.
		\end{enumerate}
		
		Moreover, if $\mathbf{V} = (V_1,\dots,V_n)$ is a $\mathbf{\Theta}_n$-isometric dilation of $\mathbf{T}$ on $\mathcal{K}$, then the operators $V_1, \dots, V_n$ are of the form given in \eqref{V_i} and \eqref{V_n}, and satisfy conditions \emph{(1)}--\emph{(4)}.
	\end{thm}
	
	The set of operators $\{A^{(i)}_l \in \mathcal{B}(\mathcal{D}_{T_n}) : 1 \le i \le n-1, 0 \le l \le p\}$ in Theorem \ref{Conditional Dilation} plays a pivotal role for constructing the conditional dilation of a $\mathbf{\Theta}_n$-contraction. As per our requirement, we rewrite the set as $\{A^{(i)}_l \in \mathcal{B}(\mathcal{D}_{T_n}) : 1 \le i \le n-1, 0 \le l \le p-1\}$ satisfying the conditions
	\begin{equation}\label{Compatibility Condition}
		\begin{aligned}
			A^{(i)*}_l = A^{(n-i)}_{p-l}, \quad 1 \le i \le n-1, \quad 0 \le l \le p-1,
		\end{aligned}
	\end{equation}
	\begin{equation}\label{Commutativity Condition}
		\begin{aligned}
			\displaystyle\sum_{l=0}^{k}[A^{(i)}_l, A^{(j)}_{k-l}] = 0, \quad 0 \le k \le 2p \quad 1 \le i, j \le n-1, \quad [A, B] = AB - BA,
		\end{aligned}
	\end{equation}
	and
	\begin{equation}\label{Condition}
		\begin{aligned}
			D_{T_n}T_i = \displaystyle\sum_{l=0}^p A^{(i)}_lD_{T_n}T^l_n, \quad 1 \le i \le n-1.
		\end{aligned}
	\end{equation}
	when we consider $(A^{(i)}_0, \dots, A^{(i)}_{p-1}), 1 \le i \le n-1$ as tuples of fundamental operators for $\mathbf{T}$.
	
	In \cite[Example 2.12, Example 3.3]{Mandal}, S. Mandal and A. Pal showed that the condition
	\begin{equation*}
		\begin{aligned}
			F_iF^*_{n-j} - F_jF^*_{n-i} = F^*_{n-j}F_i - F^*_{n-i}F_j, \quad 1 \le i, j \le n-1
		\end{aligned}
	\end{equation*}
	is not necessary, in general. M. Bhowmik in \cite[Theorem 2.3]{Bhowmik} presented a family of tetrablock contractions (see \cite{Bhattacharyya}) that has tetrablock-unitary dilation, whereas their fundamental operators are not commuting with each other. It is also shown in \cite[Theorem 2.4]{Paul 3} that there are $\Gamma_n$-contractions that extend to a $\Gamma_n$-isometry, while its fundamental operators satisfy none of the conditions assumed in \cite[Theorem 6.6, Corollary 6.67]{A. Pal}. All these developments and Theorem \ref{Conditional Dilation} motivate us to work on the following problem.
	
	\vspace{0.3cm}
	
	\noindent{\textbf{Problem 1.}}
	\textit{Can the conditions in \eqref{Compatibility Condition} and \eqref{Commutativity Condition} be waived to obtain a $\mathbf{\Theta}_n$-isometric dilation? Equivalently, are the conditions in \eqref{Compatibility Condition} and \eqref{Commutativity Condition} necessary for the existence of a $\mathbf{\Theta}_n$-isometric dilation?}
	
	\vspace{0.3cm}
	
	An open question is given in \cite{Paul 1} which we quote in the following as a problem.
	
	\vspace{0.3cm}
	
	\noindent{\textbf{Problem 2.}}
	\textit{Do the fundamental operators associated with a $\mathbf{\Theta}_n$-contraction satisfy equation \eqref{Condition}?}
	
	\vspace{0.3cm}
	
	The notion of admissible fundamental operators first introduced in \cite{Lata} for the $\Gamma$-contraction and tetrablock contraction. B. Bisai and S. Pal in \cite[Theorem 3.4]{S. Pal} generalized the work for $\Gamma_n$-contraction. We consider the following problem in the context of $\mathbf{\Theta}_n$-contraction.
	
	\vspace{0.3cm}
	
	\noindent{\textbf{Problem 3.}}
	\textit{For a given contraction $T_n$ and a set of bounded operators $\{A^{(i)}_l : 1 \le i \le n-1, 0 \le l \le p-1\}$ acting on $\mathcal{D}_{T^*_n}$ satisfying \eqref{Compatibility Condition} and \eqref{Commutativity Condition}, does there exist a $\mathbf{\Theta}_n$-contraction such that $(A^{(i)}_1, \dots, A^{(i)}_{p-1}), 1 \le i \le n-1$ are its fundamental operators?}
	
	\vspace{0.3cm}
	
	We use the technique \cite[Theorem 2.3]{Bhowmik} for answering some of the problems of this paper.
	
	The paper is outlined as follows: In Section \ref{Non-Uniqueness}, we discuss Problem 1, where we show that the conditions in \eqref{Compatibility Condition} and \eqref{Commutativity Condition} are not necessary for the existence of a $\mathbf{\Theta}_n$-isometric dilation. The section also provides a negative answer to the question in Problem 2. Section \ref{Model and Unitary Invariant} shows that the characteristic function for a $\mathbf{\Theta}_n$-contraction satisfying \eqref{Condition *} is a unitary invariant. Moreover, in Section \ref{Admissible Fundamental Operators} we study Problem 3. It demonstrates the existence of a pure $\mathbf{\Theta}_n$-contraction such that under certain conditions, $\{(A^{(i)}_1, \dots, A^{(i)}_{p-1}) : 1 \le i \le n-1\}$ are its fundamental operators.

	\section{Non-Uniqueness of Fundamental Operators and its Consequences}\label{Non-Uniqueness}
	
	The non-uniqueness of fundamental operators of a $\mathbf{\Theta}_n$-contraction is shown in \cite[Subsection 2.1]{Paul 1}. It is also described in that paper when the fundamental equations produce a unique solution (see \cite[Theorem 2.13]{Paul 1}). This section begins with another family of non-zero $\mathbf{\Theta}_n$-contractions whose fundamental equations do not have a unique solution. The absence of uniqueness of the tuples of fundamental operators is crucial in solving Problem 1 and Problem 2. Before presenting the main results of this section, we recall a result by Knese.
	
	\begin{prop}[{\cite[Proposition 1.1]{Knese}}]\label{Prop Knese}
		If $f : \mathbb{D}^n \to \mathbb{D}$ is holomorphic, then
		\begin{equation*}
			\begin{aligned}
				\sum_{i=1}^{n} (1 - |z_i|^2)\left|\frac{\partial f}{\partial z_i}(z)\right| &\le 1 - |f(z)|^2
			\end{aligned}
		\end{equation*}
		for any $z = (z_1, \dots, z_n) \in \mathbb{D}^n$.
	\end{prop}
	
	We denote the set $\{(r\theta_1, \dots, r\theta_n) : (\theta_1, \dots, \theta_n) \in \mathbf{\Theta}_n\}$ by $r\mathbb{D}^n$. Observe that $|f(z, (0, \dots, 0))| = 0 < 1$. Thus, $\mathbf{0} = (0, \dots, 0) \in \mathbf{\Theta}_n$. Then $\mathbf{0}$ is an interior point of $\mathbf{\Theta}_n$. This implies that there exists $r \in (0, 1)$ such that $r\overline{\mathbb{D}}^n$ is entirely contained in $\mathbf{\Theta}_n$. We therefore proved the following lemma.
	
	\begin{lem}\label{Lem 1}
		There exists an $r \in (0, 1)$ such that $r\overline{\mathbb{D}}^n \subseteq \mathbf{\Theta}_n$.
	\end{lem}
	
	Suppose $|\delta_j| < r^{1/m}/M$ for $1 \le j \le n-1$, where $M = \max\{k(i) : 1 \le i \le n\}$ and $\delta_n = 0$. Let $\lambda_i = \theta_i(\delta_1, \dots, \delta_n)$ for $1 \le i \le n$. Thus, $\lambda_n = 0$. We define the following nilpotent matrices:
	\begin{equation}\label{T_i, T_n}
		\begin{aligned}
			T_i &=
			\left[
			\begin{matrix}
				0 & \lambda_i\\
				0 & 0
			\end{matrix}\right], \quad 1 \le i \le n - 1, \quad
			T_n =
			\left[
			\begin{matrix}
				0 & 0\\
				0 & 0
			\end{matrix}\right].
		\end{aligned}
	\end{equation}
	
	The following theorem is one of the main results of this paper.
	
	\begin{thm}\label{Thm 1}
		Let $r$ be as in Lemma \ref{Lem 1} and $|\delta_j| < r^{1/m}/M$ for $1 \le j \le n-1$, where $M = \max\{k(i) : 1 \le i \le n\}$ and $\delta_n = 0$. Let $\mathbf{T} = (T_1, \dots, T_n)$ be as in \eqref{T_i, T_n}. Then
		\begin{enumerate}
			\item $\mathbf{T}$ is a $\mathbf{\Theta}_n$-contraction that possesses $\mathbf{\Theta}_n$-isometric dilation.
			
			\item \eqref{Fundamental} does not have a unique solution, for $m > 1, p > 1$.
			
			\item the conditions in \eqref{Compatibility Condition} and \eqref{Commutativity Condition} need not be satisfied, in general.
		\end{enumerate}
	\end{thm}
	
	\begin{proof}
		First we prove that $\mathbf{T}$ is a $\mathbf{\Theta}_n$ contraction. Note that the operators $T_1, \dots, T_n$ commute with each other in $M_2(\mathbb{C})$, as $T_iT_j = T_jT_i = 0, 1 \le i < j \le n$. As $T_1, \dots, T_{n-1}$ are nilpotent and $T_n = 0$, the only eigenvalue of $T_i$'s is $0$. By the fact that the joint spectrum of a tuple of commuting matrices is the set of joint eigenvalues of the matrices, we have that the only joint eigenvalue of the tuple $(T_1, \dots, T_n)$ is $(0, \dots, 0)$ with a joint eigenvector
		$\left[\begin{smallmatrix}
			x\\
			y
		\end{smallmatrix}\right]$ in $\mathbb{C}^2$. One can easily check that
		$\left[\begin{smallmatrix}
			1\\
			0
		\end{smallmatrix}\right]$ is one of the joint eigenvectors of $(T_1, \dots, T_n)$. Therefore, $\sigma(T_1, \dots, T_n) = \{(0, \dots, 0)\} \subseteq \mathbf{\Theta}_n$.
		
		Suppose that $f \in \mathcal{O}(\mathbf{\Theta}_n)$ such that $f(0, \dots, 0) = \alpha \in \mathbb{D} \setminus \{0\}$. Indeed, if $f(0, \dots, 0) = c \in \mathbb{D} \setminus \{0\}$, then $(\Psi_{\alpha} \circ f)(0, \dots, 0) = 0$, where $\Psi_{\alpha}(z) = e^{i\theta}\frac{z - \alpha}{1 - \overline{\alpha}z}, z \in \mathbb{D}, \theta \in \mathbb{R}$ is an automorphism of $\mathbb{D}$.
		
		As $T^2_i = 0$ for $1 \le i \le n-1$ and $T_n = 0$, by Taylor series expansion of $f$ and the holomorphic functional calculus, we obtain
		\begin{equation*}
			\begin{aligned}
				f(T_1, \dots, T_n) &=
				\left[
				\begin{matrix}
					0 & \displaystyle\sum_{i=1}^{n-1} \lambda_i \frac{\partial f}{\partial z_i}(0, \dots, 0)\\
					0 & 0
				\end{matrix}\right].
			\end{aligned}
		\end{equation*}
		
		We consider two cases. One is $||f||_{\infty, \overline{\mathbf{\Theta}}_n} = 1$ and the other is $||f||_{\infty, \overline{\mathbf{\Theta}}_n} \ne 1$.\\
		\noindent{\textbf{Case 1:}} Suppose $||f||_{\infty, \overline{\mathbf{\Theta}}_n} = 1$. Applying Proposition \ref{Prop Knese} on $f(rz_1, \dots, rz_n)$ we obtain
		\begin{equation}
			\begin{aligned}
				||f(T_1, \dots, T_n)||
				= \left|\displaystyle\sum_{i=1}^{n-1} \lambda_i \frac{\partial f}{\partial z_i}(0, \dots, 0)\right|
				&\le r\left(\displaystyle\sum_{i=1}^{n-1} \left|\frac{\partial f}{\partial z_i}(0, \dots, 0)\right|\right)\\
				&\le r\left(\displaystyle\sum_{i=1}^{n} \left|\frac{\partial f}{\partial z_i}(0, \dots, 0)\right|\right)\\
				&\le 1 - |f(0, \dots, 0)|^2\\
				&= 1 = ||f||_{\infty, \overline{\mathbf{\Theta}}_n}.
			\end{aligned}
		\end{equation}
		
		\noindent{\textbf{Case 2:}} If $||f||_{\infty, \overline{\mathbf{\Theta}}_n} \neq 1$, then performing similar technique to the function $\widetilde{f} = \frac{f}{||f||_{\infty, \overline{\mathbf{\Theta}}_n}}$ we have
		\begin{equation}\label{Main 1}
			\begin{aligned}
				||\widetilde{f}(T_1, \dots, T_n)||
				&\le 1.
			\end{aligned}
		\end{equation}
		Hence, from \eqref{Main 1} we deduce that
		\begin{equation*}
			\begin{aligned}
				||f(T_1, \dots, T_n)|| &\le ||f||_{\infty, \overline{\mathbf{\Theta}}_n}.
			\end{aligned}
		\end{equation*}
		Both cases show that $(T_1, \dots, T_n)$ is a $\mathbf{\Theta}_n$-contraction. The existence of $\mathbf{\Theta}_n$-isometric dilation follows form Theorem \ref{2 times 2 Theta_n-Contraction}. This prove $(1)$.
		
		$(2):$ In order to show that \eqref{Fundamental} does not have a unique solution, we observe that, since $T_n = 0$,
		\begin{equation*}
			\begin{aligned}
				D^2_{T_n} &= I - T^*_nT_n =
				\left[
				\begin{matrix}
					1 & 0\\
					0 & 1
				\end{matrix}\right].
			\end{aligned}
		\end{equation*}
		Thus, we have
		\begin{equation*}
			\begin{aligned}
				D_{T_n} &=
				\left[
				\begin{matrix}
					1 & 0\\
					0 & 1
				\end{matrix}\right].
			\end{aligned}
		\end{equation*}
		As $T^l_n = 0, 1 \le l \le p$, then \eqref{Fundamental} boils down to
		\begin{equation*}
			\begin{aligned}
				T_i &= D_{T_n}A^{(i)}_0D_{T_n}.
			\end{aligned}
		\end{equation*}
		One can easily deduce that $A^{(i)}_0 = T_i$. Now it follows from here that for any $p-1$ bounded operators $A^{(i)}_1, \dots, A^{(i)}_{p-1}$, the tuple of operators $(T_i, A^{(i)}_1, \dots, A^{(i)}_{p-1})$ is a solution of \eqref{Fundamental}. Therefore, for $m > 1, p > 1$, \eqref{Fundamental} does not have a unique solution.
		
		As an immediate consequence of part $(2)$, it is easy to see that the conditions in \eqref{Compatibility Condition} can only be satisfied for choosing $A^{(i)*}_l = A^{(n-i)}_{p-l}$ for $1 \le i \le n-1$ and $0 \le l \le p-1$. For proving the assertion $(3)$, we take the $(m, p, n) = (2, 2, 3)$. Due to the non-uniqueness of the fundamental operators, we assume
		\begin{equation*}
			\begin{aligned}
				A^{(1)}_1 =
				\left[
				\begin{matrix}
					1 & 0\\
					1 & 1
				\end{matrix}\right], \quad \text{and} \quad
				A^{(2)}_1 =
				\left[
				\begin{matrix}
					1 & 0\\
					0 & 1
				\end{matrix}\right].			\end{aligned}
		\end{equation*}
		Thus, $A^{(1)*}_1 \ne A^{(2)*}_1$. Since $A^{(i)}_0 = T_i =
		\left[
		\begin{matrix}
			0 & \lambda_i\\
			0 & 0
		\end{matrix}\right], 1 \le i \le 2$, then by routine calculation we obtain
		\begin{equation*}
			\begin{aligned}
				\sum_{l=0}^{1} [A^{(1)}_l, A^{(2)}_{1-l}] &=
				\left[
				\begin{matrix}
					-\lambda_2 & 0\\
					0 & \lambda_2
				\end{matrix}\right] \ne 0,
			\end{aligned}
		\end{equation*}
		proving $(3)$. This completes the proof.
	\end{proof}
	
	Next, we will discuss Problem 2 and provide a negative answer. Notice that any $\mathbf{\Theta}_n$-isometry satisfies \eqref{Condition}. Hence, we need to work with a $\mathbf{\Theta}_n$-contraction that is not a $\mathbf{\Theta}_n$-isometry. Consider $m = p = 2, n = 3$ and choose $\delta_1, \delta_2, \delta_3$ in $r^{1/3}\mathbb{D}$ such that none of them is zero and
	\begin{equation*}
		\begin{aligned}
			\delta^2_1 + \delta^2_2 + \delta^2_3 = 0 = \delta^2_1\delta^2_2 + \delta^2_2\delta^2_3 + \delta^2_3\delta^2_1.
		\end{aligned}
	\end{equation*}
	Let $z_1z_2z_3 = \lambda$. Then $(0, 0, \lambda) \in \overline{\mathbf{\Theta}}_n$. We consider the following triple of operators:
	\begin{equation}\label{T_1, T_2, T_3}
		\begin{aligned}
			(T_1, T_2, T_3) &=
			\left(
			\left[
			\begin{matrix}
				0 & 0\\
				0 & 0
			\end{matrix}\right],
			\left[
			\begin{matrix}
				0 & 0\\
				0 & 0
			\end{matrix}\right],
			\left[
			\begin{matrix}
				0 & \lambda\\
				0 & 0
			\end{matrix}\right]\right), \quad \lambda \in r\mathbb{D}.
		\end{aligned}
	\end{equation}
	
	\begin{thm}\label{Thm 2}
		Let $m = p = 2$ and $n = 3$. Suppose $r$ is as of Lemma \ref{Lem 1} and $0 < |\delta_i| < r^{1/3}$ for $1 \le i \le 3$, such that $\delta^2_1 + \delta^2_2 + \delta^2_3 = 0 = \delta^2_1\delta^2_2 + \delta^2_2\delta^2_3 + \delta^2_3\delta^2_1$ and $\delta_1\delta_2\delta_3 = \lambda$. Let $\mathbf{T} = (T_1, T_2, T_3)$ be as in \eqref{T_1, T_2, T_3}. Then
		\begin{enumerate}
			\item $\mathbf{T}$ is a $\mathbf{\Theta}_n$-contraction.
			
			\item the fundamental operators of $\mathbf{T}$ do not satisfy \eqref{Condition}.
		\end{enumerate}
	\end{thm}
	
	\begin{proof}
		We only prove $(2)$ as the proof of $(1)$ is similar to that of part $(1)$ of Theorem \ref{Thm 1}. For proving $(2)$, note that
		\begin{equation*}
			\begin{aligned}
				D_{T_3} &=
				\left[
				\begin{matrix}
					1 & 0\\
					0 & \sqrt{1 - |\lambda|^2}
				\end{matrix}\right].
			\end{aligned}
		\end{equation*}
		We now find tuples of fundamental operators of $\mathbf{T}$. Observe that
		\begin{equation*}
			\begin{aligned}
				0 = T_1 - T^*_2T^2_3
				&= T^*_3D_{T_3}A^{(1)}_0D_{T_3}T_3 + D_{T_3}A^{(1)}_1D_{T_3}T_3 + D_{T_3}A^{(1)}_0D_{T_3}.
			\end{aligned}
		\end{equation*}
		Suppose that
		$A^{(1)}_0 =
		\left[
		\begin{matrix}
			x_{11} & x_{12}\\
			x_{21} & x_{22}
		\end{matrix}\right]$ and $A^{(1)}_1 =
		\left[
		\begin{matrix}
		y_{11} & y_{12}\\
		y_{21} & y_{22}
		\end{matrix}\right]$. After rigorous computation, it results that
		\begin{equation*}
			\begin{aligned}
				\left[
				\begin{matrix}
					x_{11} & \lambda y_{11} + x_{12}\sqrt{1 - |\lambda|^2}\\
					x_{21}\sqrt{1 - |\lambda|^2} & \lambda y_{21}\sqrt{1 - |\lambda|^2} + x_{22}(1 - |\lambda|^2)
				\end{matrix}\right] &=
				\left[
				\begin{matrix}
					0 & 0\\
					0 & 0
				\end{matrix}\right].
			\end{aligned}
		\end{equation*}
		Thus, we have $x_{11} = x_{21} = 0$, and
		\begin{equation}\label{Main 2}
			\begin{aligned}
				&\lambda y_{11} + x_{12}\sqrt{1 - |\lambda|^2} = 0, \quad \text{and} \quad \lambda y_{21} + x_{22}\sqrt{1 - |\lambda|^2}  = 0.
			\end{aligned}
		\end{equation}
		This implies that the solution sets of the equations in \eqref{Main 2} are
		\begin{equation*}
			\begin{aligned}
				\{\alpha(- \lambda, \sqrt{1 - |\lambda|^2}) : \alpha \in \mathbb{C}\}, \quad \text{and} \quad \{\beta(- \lambda, \sqrt{1 - |\lambda|^2}) : \beta \in \mathbb{C}\}
			\end{aligned}
		\end{equation*}
		respectively. It then follows from here that
		\begin{equation*}
			\begin{aligned}
				A^{(1)}_0 =
				\left[
				\begin{matrix}
					0 & - \alpha\lambda\\
					0 & - \beta\lambda
				\end{matrix}\right], \quad \text{and} \quad
				A^{(1)}_1 =
				\left[
				\begin{matrix}
					\alpha\sqrt{1 - |\lambda|^2} & *\\
					\beta\sqrt{1 - |\lambda|^2} & *
				\end{matrix}\right],
			\end{aligned}
		\end{equation*}
		which yields that
		\begin{equation*}
			\begin{aligned}
				D_{T_3}T_1 &=
				\left[
				\begin{matrix}
					0 & 0\\
					0 & 0
				\end{matrix}\right] \ne
				\left[
				\begin{matrix}
					0 & - \alpha\lambda + \alpha\lambda\sqrt{1 - |\lambda|^2}\\
					0 & - \beta\lambda + \beta\sqrt{1 - |\lambda|^2}
				\end{matrix}\right] =
				A^{(1)}_0D_{T_3} + A^{(1)}_1D_{T_3}T_3 + A^{(1)}_1D_{T_3}T^2_3.
			\end{aligned}
		\end{equation*}
		Therefore the tuples of fundamental operators do not satisfy \eqref{Condition}. This finishes the proof.
	\end{proof}

	\section{Fundamental Operators and Functional Model for Pure $\mathbf{\Theta}_n$-Contraction}\label{Model and Unitary Invariant}
	
	We begin our discussion by recalling some developments in one variable settings. Let $T$ and $\widetilde{T}$ be two contractions on two different Hilbert space $\mathcal{H}$ and $\widetilde{\mathcal{H}}$ respectively. The characteristic functions of $T$ and $\widetilde{T}$ are said to coincide if there exist unitary operators $U : \mathcal{D}_T \to \mathcal{D}_{\widetilde{T}}$ and $U_* : \mathcal{D}_{T^*} \to \mathcal{D}_{\widetilde{T}}$ such that the diagram
	\begin{equation}\label{Commutative Diagram}
		\begin{aligned}
			\begin{tikzcd}
				\mathcal{D}_T \arrow{r}{\Theta_T(z)} \arrow[swap]{d}{U} & \mathcal{D}_{T^*} \arrow{d}{U_*} \\
				\mathcal{D}_{\widetilde{T}} \arrow{r}{\Theta_{\widetilde{T}}(z)}& \mathcal{D}_{\widetilde{T}^*}
			\end{tikzcd}
		\end{aligned}
	\end{equation}
	commute for all $z \in \mathbb{D}$.
	
	\begin{prop}[Nagy-Foias]\label{NF Prop 11}
		Characteristic functions of unitarily equivalent contractions coincide.
	\end{prop}
	
	The converse of Proposition \ref{NF Prop 11} holds for pure contractions. We recall it in the following proposition.
	
	\begin{prop}\label{Prop 12}
		Let $T$ and $\widetilde{T}$ be two pure contractions on Hilbert spaces $\mathcal{H}$ and $\widetilde{\mathcal{H}}$ respectively. If their characteristic functions coincide then $T$ and $\widetilde{T}$ are unitarily equivalent.
	\end{prop}
	
	Sz.-Nagy and Foias in \cite{Nagy} showed that, the characteristic function of a completely non-unitary contraction is a complete unitary invariant. We quote the result in the following theorem.
	
	\begin{thm}[Nagy-Foias]\label{NF Thm 12}
		Two completely non-unitary contractions are unitarily equivalent if and only if their characteristic functions coincide.
	\end{thm}
	
	In this section, we show that for a given pure
	$\mathbf{\Theta}_n$-contraction $\mathbf{T} = (T_1, \dots, T_n)$, the characteristic function of $T_n$ is unitary invariant. Define
	$W : \mathcal{H} \to
	H^2(\mathbb{D}) \otimes \mathcal{D}_{T^*_n}$
	by
	\begin{equation}\label{W}
		\begin{aligned}
			Wh &=
			\sum_{k \ge 0} z^k \otimes D_{T^*_n}T^{*k}_nh,
			\quad h \in \mathcal{H}.
		\end{aligned}
	\end{equation}
	Since $T_n$ is a pure contraction, the standard defect-space calculation shows that $W$ is an isometry. Indeed, due to the pureness of $T_n$ we have
	\begin{equation*}
		\begin{aligned}
			||Wh||^2
			&= \sum_{k \ge 0}
			||D_{T^*_n}T^{*k}_nh||^2 = \sum_{k \ge 0}
			\left(
			||T^{*k}_nh||^2 - ||T^{*(k+1)}_nh||^2
			\right) = ||h||^2.
		\end{aligned}
	\end{equation*}
	One can show that
	\begin{equation}\label{W*}
		\begin{aligned}
			\begin{aligned}
				W^*(z^k \otimes \eta)  &=
				T^k_nD_{T^*_n}\eta,
				\quad \eta \in \mathcal{D}_{T^*_n}, k \ge 0.
			\end{aligned}
		\end{aligned}
	\end{equation}
	
	We state the following lemma. The proof can be found in \cite[Lemma 4.5]{S. Pal 3}.
	
	\begin{lem}\label{W isometry}
		Let $T_n$ be contraction. Then
		\begin{equation}\label{W Property}
			\begin{aligned}
				WW^* + M_{\Theta_{T_n}}M^*_{\Theta_{T_n}} = I_{H^2(\mathbb{D}) \otimes \mathcal{D}_{T^*_n}}.
			\end{aligned}
		\end{equation}
	\end{lem}
	
	In the following theorem, we recall the functional model for a class of pure $\mathbf{\Theta}_n$-contraction that satisfies the condition
	\begin{equation}\label{Condition *}
		\begin{aligned}
			D_{T^*_n}T^*_i &= \sum_{l=0}^{p} A^{(i)^*}_lD_{T^*_n}T^{*l}_n
		\end{aligned}
	\end{equation}
	is described \cite[Theorem 3.6]{Paul 2}.
	
	\begin{thm}\label{Functional Model}
		Let $\mathbf{T} = (T_1, \dots, T_n)$ be a pure $\mathbf{\Theta}_n$-contraction on a Hilbert space $\mathcal{H}$. Suppose that $\{A^{(i)}_l : 0 \le l \le p, 1 \le i \le n-1\}$ is a set of operators on $\mathcal{D}_{T_n}$ that satisfy \eqref{Compatibility Condition} and \eqref{Condition *}. Then
		\begin{enumerate}
			\item $T_i \cong \mathbb{P}_{\mathcal{H}_{T_n}}\left(\displaystyle\sum_{l=0}^{p} M^l_z \otimes A^{(i)}_l\right)\Bigg|_{\mathcal{H}_{T_n}}, 1 \le i \le n-1$,
			
			\item $T_n \cong \mathbb{P}_{\mathcal{H}_{T_n}}(M_z \otimes I_{\mathcal{D}_{T^*_n}})|_{\mathcal{H}_{T_n}}$,
		\end{enumerate}
		where the notation ``$\cong$'' stands for unitary equivalence, $\mathcal{H}_{T_n} =
		(H^2(\mathbb{D}) \otimes \mathcal{D}_{T^*_n}) \ominus M_{\Theta_{T_n}}(H^2(\mathbb{D}) \otimes \mathcal{D}_{T_n})$ and $\mathbb{P}_{\mathcal{H}_{T_n}}$ is the projection onto the space $\mathcal{H}_{T_n}$.
	\end{thm}
	
	As an application of Theorem \ref{Functional Model}, we derive the functional model for pure $\mathbf{\Theta}_n$-isometry (see \cite[Theorem 3.5]{Biswas 2}).
	
	\begin{cor}\label{Cor 1}
		Let $\mathbf{T} = (T_1, \dots, T_n)$ be a pure $\mathbf{\Theta}_n$-isometry on a Hilbert space $\mathcal{H}$. Then $\mathbf{T}$ is unitarily equivalent to
		$(M_{\Phi_1}, \dots, M_{\Phi_{n-1}}, M_z)$ on $H^2(\mathcal{D}_{T^*_n})$, where $M_{\Phi_i} = \displaystyle\sum_{l=0}^{p} M^l_z \otimes A^{(i)}_l$ for some $A^{(i)}_l \in \mathcal{D}_{T^*_n}$ satisfying \eqref{Compatibility Condition} and \eqref{Commutativity Condition}.
	\end{cor}
	
	\begin{proof}
		Given $\mathbf{T}$ is a pure $\mathbf{\Theta}_n$-isometry. Since $T_n$ is an isometry, then $D_{T_n} = 0$. Consequently, $\Theta_{T_n}$ is the zero operator. Thus, the model space $\mathcal{H}_{T_n}$ for $\mathbf{T}$ is identically $H^2(\mathbb{D}) \otimes \mathcal{D}_{T^*_n}$. Therefore, Theorem \ref{Functional Model} implies that
		\begin{equation*}
			\begin{aligned}
				\mathbf{T} \cong \left(\displaystyle\sum_{l=0}^{p} M^l_z \otimes A^{(1)}_l, \dots, \displaystyle\sum_{l=0}^{p} M^l_z \otimes A^{(n-1)}_l, M_z \otimes I_{\mathcal{D}_{T^*_n}}\right),
			\end{aligned}
		\end{equation*}
		where $A^{(i)}_l$'s satisfy \eqref{Compatibility Condition} and . As we can identify $H^2(\mathbb{D}) \otimes \mathcal{D}_{T^*_n}$ by $H^2(\mathcal{D}_{T^*_n})$, it is clear that $\mathbf{T} \cong (M_{\Phi_1}, \dots, M_{\Phi_{n-1}}, M_z)$, where $\Phi_i(z) = \displaystyle\sum_{l=0}^{p} A^{(i)}_lz^l$ for $1 \le i \le n-1$, with $A^{(i)}_l$'s obeying \eqref{Compatibility Condition}. Since $T_i$'s commute with each other then so are $M_{\Phi_i}$'s for $1 \le i \le n-1$, which establishes \eqref{Commutativity Condition}. This completes the proof.
	\end{proof}
	
	Notice that any $\mathbf{\Theta}_n$-isometry satisfies \eqref{Condition}. As a result, it can be easily seen that we have deduced the functional model for a pure $\mathbf{\Theta}_n$-isometry (Corollary \ref{Cor 1}) from a functional model for a pure $\mathbf{\Theta}_n$-contraction satisfying \eqref{Condition}. On the other hand, Corollary \ref{Cor 1} implies that the Hilbert space $\mathcal{E}$ in \cite[Theorem 3.5]{Biswas 2} can be taken as $\mathcal{D}_{T^*_n}$.
	
	The remaining part of this section talks about some unitary invariants related to $\mathbf{\Theta}_n$-contraction. For any operator $T$ and a unitary $U$, we call $UTU^*$ the conjugation of $T$ by $U$.
	
	\begin{thm}\label{Conjugation Fundamental}
		Let $\mathbf{T} = (T_1, \dots, T_n)$ and $\widetilde{\mathbf{T}} = (\widetilde{T}_1, \dots, \widetilde{T}_n)$ be two $\mathbf{\Theta}_n$-contraction on Hilbert spaces $\mathcal{H}$ and $\widetilde{\mathcal{H}}$, respectively, such that $\mathbf{T}$ and $\widetilde{\mathbf{T}}$ are unitarily equivalent. If $(A^{(i)}_0, \dots, A^{(i)}_{p-1})$ is a tuple of fundamental operators of $\mathbf{T}$, then the conjugation of $(A^{(i)}_0, \dots, A^{(i)}_{p-1})$ by the same unitary involved in the unitary equivalence of $\mathbf{T}$ and $\widetilde{\mathbf{T}}$ are the tuples of fundamental operators of $\widetilde{\mathbf{T}}$. Moreover, the characteristic functions of $T_n$ and $\widetilde{T}_n$ coincide.
	\end{thm}
	
	\begin{proof}
		Given $\mathbf{T}$ and $\widetilde{\mathbf{T}}$ are unitarily equivalent. Then there exists unitary $U : \mathcal{H} \to \widetilde{\mathcal{H}}$ such that
		\begin{equation}\label{UI 1}
			\begin{aligned}
				UT_i = \widetilde{T}_iU, \quad 1 \le i \le n-1, \quad UT_n = \widetilde{T}_nU.
			\end{aligned}
		\end{equation}
		From here it implies that
		\begin{equation}\label{UI 2}
			\begin{aligned}
				UD^2_{T_n} = U(I - T^*_nT_n) = (I - \widetilde{T}^*_n\widetilde{T}_n)U = D^2_{\widetilde{T}_n}U, \quad i.e., \quad UD_{T_n} = D_{\widetilde{T}_n}U.
			\end{aligned}
		\end{equation}
		We consider $\widetilde{U} = U|_{\mathcal{D}_{T_n}} : \mathcal{D}_{T_n} \to \mathcal{D}_{\widetilde{T}_n}$. Thus $\widetilde{U} \in \mathcal{B}(\mathcal{D}_{T_n}, \mathcal{D}_{\widetilde{T}_n})$. Then $\widetilde{U}$ is a unitary.
		
		Fix any $i$ in $\{1 , \dots, n-1\}$. Since $(A^{(i)}_0, \dots, A^{(i)}_{p-1})$ a tuple of fundamental operators of \eqref{Fundamental}, then we have
		\begin{equation}\label{UI 3}
			\begin{aligned}
				T_i - T^*_{n-i}T^p_n
				&= \sum_{l=0}^{p-1} T^{*p-1-l}_n D_{T_n}A^{(i)}_lD_{T_n} T^{p-1}_n + \sum_{l=0}^{p-2} T^{*p-2-l}_n D_{T_n}A^{(i)}_lD_{T_n} T^{p-2}_n + \dots\\
				& \quad + \sum_{l=0}^{1} T^{*1-l}_n D_{T_n}A^{(i)}_lD_{T_n} T_n + D_{T_n}A^{(i)}_0D_{T_n}.
			\end{aligned}
		\end{equation}
		Conjugating \eqref{UI 3} both side by $\widetilde{U}$ we get
		\begin{equation*}
			\begin{aligned}
				\widetilde{U}(T_i - T^*_{n-i}T^p_n)\widetilde{U}^* &= \sum_{l=0}^{p-1} \widetilde{U}T^{p-1-l}_n D_{T^*_n}A^{(i)}_lD_{T^*_n}T^{*p-1}_n\widetilde{U}^* + \sum_{l=0}^{p-2} \widetilde{U}T^{p-2-l}_n D_{T^*_n}A^{(i)}_lD_{T^*_n}T^{*p-2}_n\widetilde{U}^* + \dots\\
				& \quad + \sum_{l=0}^{1} \widetilde{U}T^{1-l}_n D_{T^*_n}A^{(i)}_lD_{T^*_n}T^*_n\widetilde{U}^* + \widetilde{U}D_{T^*_n}A^{(i)}_0D_{T^*_n}\widetilde{U}^*.
			\end{aligned}
		\end{equation*}
		Using \eqref{UI 1} and \eqref{UI 2} we derive
		\begin{equation*}
			\begin{aligned}
				\widetilde{T}_i - \widetilde{T}^*_{n-i}\widetilde{T}^p_n
				&= \sum_{l=0}^{p-1} \widetilde{T}^{*p-1-l}_n D_{\widetilde{T}_n}(\widetilde{U}A^{(i)}_l\widetilde{U}^*)D_{\widetilde{T}_n} \widetilde{T}^{p-1}_n + \sum_{l=0}^{p-2} \widetilde{T}^{*p-2-l}_n D_{\widetilde{T}_n}(\widetilde{U}A^{(i)}_l\widetilde{U}^*)D_{\widetilde{T}_n} \widetilde{T}^{p-2}_n + \dots\\
				& \quad + \sum_{l=0}^{1} \widetilde{T}^{*1-l}_n D_{\widetilde{T}_n}(\widetilde{U}A^{(i)}_l\widetilde{U}^*)D_{\widetilde{T}_n} T_n + D_{\widetilde{T}_n}(\widetilde{U}A^{(i)}_l\widetilde{U}^*)D_{\widetilde{T}_n}.
			\end{aligned}
		\end{equation*}
		This implies that $(\widetilde{U}A^{(i)}_0\widetilde{U}^*, \dots, \widetilde{U}A^{(i)}_{p-1}\widetilde{U}^*)$ is a tuple of fundamental operators of $\widetilde{\mathbf{T}}$. Since $\mathbf{T}$ is unitarily equivalent to $\widetilde{\mathbf{T}}$, then so are $T_n$ and $\widetilde{T}_n$. Therefore, by Proposition \ref{NF Prop 11}, the characteristic functions of $T_n$ and $\widetilde{T}_n$ coincide. This completes the proof.
	\end{proof}
	
	\begin{rem}
		Although in Theorem \ref{Conjugation Fundamental} we proved that the unitary conjugation of a tuple of fundamental operators of a $\mathbf{\Theta}_n$-contraction $\mathbf{T}$ is a tuple of fundamental operators of the $\mathbf{\Theta}_n$-contraction $\widetilde{\mathbf{T}}$ to which $\mathbf{T}$ is unitarily equivalent, but due to the unavailability of the uniqueness of the fundamental operators, we cannot conclude that the tuples of fundamental operators of two unitarily equivalent $\mathbf{\Theta}_n$-contractions are also unitarily equivalent.
	\end{rem}
	
	The following theorem establishes the unitary equivalence of two pure $\mathbf{\Theta}_n$-contractions $\mathbf{T}$ and $\widetilde{\mathbf{T}}$ satisfying \eqref{Condition *} in the presence of the coincidence of the characteristic functions of $T_n$ and $\widetilde{T}_n$ and the unitary equivalence of the tuples of fundamental operators.
	
	\begin{thm}\label{Unitary Invariants}
		Let $\mathbf{T} = (T_1, \dots, T_n)$ and $\widetilde{\mathbf{T}} = (\widetilde{T}_1, \dots, \widetilde{T}_n)$ be two pure $\mathbf{\Theta}_n$-contraction on the Hilbert spaces $\mathcal{H}$ and $\widetilde{\mathcal{H}}$, respectively, satisfying \eqref{Condition *}. Suppose that the characteristic functions of $T_n$ and $\widetilde{T}_n$ coincide and the tuple of fundamental operators $(A^{(i)}_0, \dots, A^{(i)}_{p-1}), 1 \le i \le n-1$ of $\mathbf{T}$ is unitarily equivalent to the tuple of fundamental operators $(\widetilde{A}^{(i)}_0, \dots, \widetilde{A}^{(i)}_{p-1}), 1 \le i \le n-1$ of $\widetilde{\mathbf{T}}$ by the same unitary involved in the coincidence of the characteristic functions of $T_n$ and $\widetilde{T}_n$. Then $\mathbf{T}$ and $\widetilde{\mathbf{T}}$ are unitarily equivalent.
	\end{thm}
	
	\begin{proof}
		Let $U : \mathcal{D}_{T_n} \to \mathcal{D}_{\widetilde{T}_n}$ and $U_* : \mathcal{D}_{T^*_n} \to \mathcal{D}_{\widetilde{T}^*_n}$ be unitaries such that for $ 1 \le i \le n-1$ and $0 \le l \le p-1$,
		\begin{equation}\label{UE 1}
			\begin{aligned}
				&U_*A^{(i)}_l = \widetilde{A}^{(i)}_lU_*, \quad \text{and} \quad U_*\Theta_{T_n}(z) = \Theta_{\widetilde{T}_n}(z)U, \quad z \in \mathbb{D}.
			\end{aligned}
		\end{equation}
		We define $\widetilde{U}_* : H^2(\mathbb{D}) \otimes \mathcal{D}_{T^*_n} \to H^2(\mathbb{D}) \otimes \mathcal{D}_{\widetilde{T}^*_n}$ by
		\begin{equation*}
			\begin{aligned}
				&\widetilde{U}_* = \id_{H^2(\mathbb{D})} \otimes U_*, \quad i.e., \quad \widetilde{U}_*(z^k \otimes \eta) = z^k \otimes U_*\eta, \quad k \in \mathbb{N} \cup \{0\}, \eta \in \mathcal{D}_{T^*_n}.
			\end{aligned}
		\end{equation*}
		Note that for all $f \in H^2(\mathbb{D}) \otimes \mathcal{D}_{T^*_n}$, and $z \in \mathbb{D}$ we have
		\begin{equation}\label{UE 2}
			\begin{aligned}
				\widetilde{U}_*(M_{\Theta_{T_n}}f(z)) &= \widetilde{U}_*(\Theta_{T_n}(z)f(z))
				= U_*(\Theta_{T_n}(z)f(z)
				= \Theta_{\widetilde{T}_n}(z)Uf(z) = M_{\Theta_{\widetilde{T}_n}}(Uf(z)).
			\end{aligned}
		\end{equation}
		Thus, from \eqref{UE 2} we get $\widetilde{U}_*(\Ran M_{\Theta_{T_n}}) = \Ran M_{\Theta_{\widetilde{T}_n}}$; i.e., $\widetilde{U}_*$ maps $\Ran M_{\Theta_{T_n}}$ onto $\Ran M_{\Theta_{\widetilde{T}_n}}$. Since $\widetilde{U}_*$ is unitary, then
		\begin{equation}\label{UE 3}
			\begin{aligned}
				\widetilde{U}_*(\mathcal{H}_{T_n}) &= \widetilde{U}_*((H^2(\mathbb{D}) \otimes \mathcal{D}_{T^*_n}) \ominus M_{\Theta_{T_n}}(H^2(\mathbb{D}) \otimes \mathcal{D}_{T_n}))\\
				&= \widetilde{U}_*((\Ran M_{\Theta_{T_n}})^{\perp})\\
				&= (\widetilde{U}_*(\Ran M_{\Theta_{T_n}}))^{\perp}\\
				&= (\Ran M_{\Theta_{\widetilde{T}_n}}))^{\perp}\\
				&= (H^2(\mathbb{D}) \otimes \mathcal{D}_{\widetilde{T}^*_n}) \ominus M_{\Theta_{\widetilde{T}_n}}(H^2(\mathbb{D}) \otimes \mathcal{D}_{\widetilde{T}_n})\\
				&= \mathcal{H}_{\widetilde{T}_n}.
			\end{aligned}
		\end{equation}
		Hence, \eqref{UE 3} implies that $\mathcal{H}_{\widetilde{T}_n} = \widetilde{U}_*(\mathcal{H}_{T_n})$. 
		
		For showing $\mathbf{T} \cong \widetilde{\mathbf{T}}$, it is enough to show the following:
		\begin{enumerate}
			\item $\widetilde{U}_*\left(\displaystyle\sum_{l=0}^{p} M^l_z \otimes A^{(i)}_l\right)^* = \left(\displaystyle\sum_{l=0}^{p} M^l_z \otimes \widetilde{A}^{(i)}_l\right)^*\widetilde{U}_*$, for $1 \le i \le n-1$,
			
			\item $\widetilde{U}_*(M_z \otimes I_{\mathcal{D}_{T^*_n}}) = (M_z \otimes I_{\mathcal{D}_{\widetilde{T}^*_n}})\tilde{U}_*$.
		\end{enumerate}
		Now observe that
		\begin{equation}\label{UE 4}
			\begin{aligned}
				\widetilde{U}_*\left(\sum_{l=0}^{p} M^l_z \otimes A^{(i)}_l\right)^*
				&= (\id_{H^2(\mathbb{D})} \otimes U_*)\left(\sum_{l=0}^{p} M^{*l}_z \otimes A^{(i)*}_l\right)\\
				&= \sum_{l=0}^{p} M^{*l}_z \otimes U_*A^{(i)}_l\\
				&= \sum_{l=0}^{p} M^l_z \otimes \widetilde{A}^{(i)*}_lU_*\\
				&= \left(\sum_{l=0}^{p} M^{*l}_z \otimes \widetilde{A}^{(i)*}_l\right)(\id_{H^2(\mathbb{D})} \otimes U_*)\\
				&= \left(\sum_{l=0}^{p} M^l_z \otimes \widetilde{A}^{(i)}_l\right)^*\widetilde{U}_*, \quad 1 \le i \le n-1.
			\end{aligned}
		\end{equation}
		On the other hand, $(M_z \otimes I_{\mathcal{D}_{T^*_n}})$ and $(M_z \otimes I_{\mathcal{D}_{\widetilde{T}^*_n}})$ are unitarily equivalent. Indeed, note that
		\begin{equation}\label{UE 5}
			\begin{aligned}
				\widetilde{U}_*(M_z \otimes I_{\mathcal{D}_{T^*_n}})
				&= (\id_{H^2(\mathbb{D})} \otimes U_*)(M_z \otimes I_{\mathcal{D}_{T^*_n}})\\
				&= (M_z \otimes I_{\mathcal{D}_{\widetilde{T}^*_n}})(\id_{H^2(\mathbb{D})} \otimes U_*)\\
				&= (M_z \otimes I_{\mathcal{D}_{\widetilde{T}^*_n}})\tilde{U}_*.
			\end{aligned}
		\end{equation}
		It then follows from \eqref{UE 4} and \eqref{UE 5} that $\mathcal{H}_{\widetilde{T}_n}$ is joint co-invariant subspace of $\left(\displaystyle\sum_{l=0}^{p} M^l_z \otimes \widetilde{A}^{(i)}_l\right)$ for $1 \le i \le n-1$ and $(M_z \otimes I_{\mathcal{D}_{\widetilde{T}^*_n}})$ respectively. Therefore
		\begin{equation*}
			\begin{aligned}
				\mathbb{P}_{\mathcal{H}_{T_n}}\left(\sum_{l=0}^{p} M^l_z \otimes A^{(i)}_l\right)\Bigg|_{\mathcal{H}_{T_n}} \quad \cong \quad \mathbb{P}_{\mathcal{H}_{\widetilde{T}_n}}\left(\sum_{l=0}^{p} M^l_z \otimes \widetilde{A}^{(i)}_l\right)\Bigg|_{\mathcal{H}_{\widetilde{T}_n}}, \quad 1 \le i \le n-1,
			\end{aligned}
		\end{equation*}
		and
		\begin{equation*}
			\begin{aligned}
				\mathbb{P}_{\mathcal{H}_{T_n}}(M_z \otimes I_{\mathcal{D}_{T^*_n}})|_{\mathcal{H}_{T_n}} \quad \cong \quad
				\mathbb{P}_{\mathcal{H}_{\widetilde{T}_n}}(M_z \otimes I_{\mathcal{D}_{\widetilde{T}^*_n}}),
			\end{aligned}
		\end{equation*}
		By the functional model for a class of pure $\mathbf{\Theta}_n$-contraction satisfying \eqref{Condition *} we deduce that $\mathbf{T} \cong \widetilde{\mathbf{T}}$. This completes the proof.
	\end{proof}
	
	While using Theorem \ref{Conjugation Fundamental} and Theorem \ref{Unitary Invariants}, we have shown that the characteristic function is an invariant for a class of pure $\mathbf{\Theta}_n$-contraction satisfying \eqref{Condition *}, a complete set of unitary invariants for that class of pure $\mathbf{\Theta}_n$-contraction is still unknown. We therefore have the following open question.
	
	\vspace{0.3cm}
	
	\noindent{\textbf{Open Question.}} What is the complete set of unitary invariant for the class of pure $\mathbf{\Theta}_n$-contraction satisfying \eqref{Condition *}?

	\section{Admissible Fundamental Operators of $\mathbf{\Theta}_n$-Contraction}\label{Admissible Fundamental Operators}
	
	This section examines Problem 3. We produce a partial solution of this problem whenever $T_n$ is a pure contraction. For this purpose, we consider the following fundamental equation for the $n$-tuple $\mathbf{T}^* = (T^*_1, \dots, T^*_n)$:
	\begin{equation}\label{Fundamental *}
		\begin{aligned}
			T^*_i - T_{n-i}T^{*p}_n
			&= \sum_{l=0}^{p-1} T^{p-1-l}_n D_{T^*_n}A^{(i)}_lD_{T^*_n}T^{*p-1}_n + \sum_{l=0}^{p-2} T^{p-2-l}_n D_{T^*_n}A^{(i)}_lD_{T^*_n}T^{*p-2}_n + \dots\\
			& \quad + \sum_{l=0}^{1} T^{1-l}_n D_{T^*_n}A^{(i)}_lD_{T^*_n}T^*_n + D_{T^*_n}A^{(i)}_0D_{T^*_n}, \quad 1 \le i \le n-1.
		\end{aligned}
	\end{equation}
	
	We now provide the main result of this section.
	
	\begin{thm}\label{Admissible}
		Let $T_n$ be a pure contraction on a Hilbert space $\mathcal{H}$, and $\{\widetilde{A}^{(i)}_l : 1 \le i \le n-1, 0 \le l \le p-1\}$ acting $\mathcal{D}_{T^*_n}$ that satisfy \eqref{Compatibility Condition} and \eqref{Commutativity Condition}. Let
		\begin{equation}\label{Gamma_n-1}
			\begin{aligned}
				\left(\gamma_1\left(\sum_{l=0}^{p} \widetilde{A}^{(1)*}_lz^l\right), \dots, \gamma_{n-1}\left(\sum_{l=0}^{p} \widetilde{A}^{(n-1)*}_lz^l\right)\right), \quad z \in \mathbb{T}
			\end{aligned}
		\end{equation}
		is a $\Gamma_{n-1}$-contraction. Suppose there exists set of bounded operators $\{A^{(i)}_l : 1 \le i \le n-1, 0 \le l \le p-1\}$ acting $\mathcal{D}_{T_n}$ such that
		\begin{equation}\label{Admissible Condition}
			\begin{aligned}
				\left(\sum_{l=0}^{p} \widetilde{A}^{(i)*}_lz^l\right)\Theta_{T_n}(z) &= \Theta_{T_n}(z)\left(\sum_{l=0}^{p}A^{(i)}_lz^l\right), \quad z \in \mathbb{D}, \quad 1 \le i \le n-1.
			\end{aligned}
		\end{equation}
		Then there exists a $\mathbf{\Theta}_n$-contraction $\mathbf{T} = (T_1, \dots, T_n)$ such that $(\widetilde{A}^{(i)}_0, \dots, \widetilde{A}^{(i)}_{p-1})$ is a tuple of fundamental operators of \eqref{Fundamental *}. In other words, $(\widetilde{A}^{(i)}_0, \dots. \widetilde{A}^{(i)}_{p-1})$ is a tuple of fundamental operator of $\mathbf{T}^* = (T^*_1, \dots, T^*_n)$.
	\end{thm}
	
	\begin{proof}
		For the given pure contraction $T_n$, let the isometry $W$ be defined as in \eqref{W}. Observe that
		\begin{equation}\label{Ad 1}
			\begin{aligned}
				M^*_zWh &= M^*_z\left(\sum_{k \ge 0} z^kD_{T^*_n}T^{*k}_nh\right) = \sum_{k \ge 0} z^kD_{T^*_n}T^{*k+1}_nh = WT^*_nh, \quad h \in \mathcal{H}.
			\end{aligned}
		\end{equation}
		This shows that $M^*_zW = WM^*_z$. Let $\Phi_i(z) = \displaystyle\sum_{l=0}^{p} \widetilde{A}^{(i)*}_lz^l$. We define
		\begin{equation}\label{Ad 2}
			\begin{aligned}
				T_i& = W^*M_{\Phi_i}W, \quad 1 \le i \le n-1.
			\end{aligned}
		\end{equation}
		Since $T_n$ is pure, from \eqref{W isometry} we obtain $(\Ran W)^{\perp} = \Ran M_{\Theta_{T_n}}$. From \eqref{Admissible Condition}, it is clear that $\Ran M_{\Theta_{T_n}}$ is invariant under $M_{\Phi_i}$. It means that $\Ran W$ is co-invariant under $M_{\Phi_i}$. In other words, $\Ran W$ is invariant under $M^*_{\Phi_i}$.
		
		In order to prove that $\mathbf{T} = (T_1, \dots, T_n)$ is a $\mathbf{\Theta}_n$-contraction, we first need to verify the commutativity of $T_i$'s. Notice that $M_{\Phi_i}$'s commute with $M_z$. This implies that for $1 \le i \le n-1$ we have $T_iT_n = T_nT_i$. Thus, we now show that $T_iT_j = T_jT_i$ for $1 \le i < j < n-1$. Note that by \eqref{Gamma_n-1}, we have that $\Phi_i$'s commute with each other. Hence, we deduce that
		\begin{equation}\label{Ad 3}
			\begin{aligned}
				T^*_iT^*_j
				&= W^*M^*_{\Phi_i}WW^*M^*_{\Phi_j}W\\
				&= W^*M^*_{\Phi_i}M^*_{\Phi_i}W, \quad (\text{since} \,\, \Ran W \,\, \text{is invariiant under} \,\, M^*_{\Phi_j})\\
				&= W^*M^*_{\Phi_j}M^*_{\Phi_i}W, \quad (\text{as} \,\, \Phi_i\Phi_j = \Phi_j\Phi_i)\\
				&= W^*M^*_{\Phi_j}WW^*M^*_{\Phi_i}W, \quad (\text{since} \,\, \Ran W \,\, \text{is invariiant under} \,\, M^*_{\Phi_i})\\
				&= T^*_jT^*_i.
			\end{aligned}
		\end{equation}
		This yields that $T_iT_j = T_jT_i$ for $1 \le i < j \le n-1$. We therefore conclude that $\mathbf{T} = (T_1, \dots, T_n)$ is a commutative operator tuple.
		
		Next we prove that $\mathbf{T}$ is a $\mathbf{\Theta}_n$-contraction. It is immediate from \eqref{Compatibility Condition} and \eqref{Gamma_n-1} that the operator tuple $(\widehat{M}_{\Phi_1}, \dots, \widehat{M}_{\Phi_{n-1}}, \widehat{M}_z)$, where $\widehat{M}_z$ is the unitary operator acting on $L^2(\mathcal{D}_{T^*_n})$, is a $\mathbf{\Theta}_n$-unitary. Thus, being the restriction of $(\widehat{M}_{\Phi_1}, \dots, \widehat{M}_{\Phi_{n-1}}, \widehat{M}_z)$ to the joint invariant subspace $H^2(\mathcal{D}_{T^*_n})$, the multiplication operator tuple $(M_{\Phi_1}, \dots, M_{\Phi_{n-1}}, M_z)$ is a $\mathbf{\Theta}_n$-isometry. Now we see that
		\begin{equation}\label{Ad 4}
			\begin{aligned}
				p(T_1, \dots, T_{n-1}, T_n) &= W^*p(M_{\Phi_1}, \dots, M_{\Phi_{n-1}}, M_z)W,
			\end{aligned}
		\end{equation}
		for any $p \in \mathbb{C}[z_1, \dots, z_n]$. As $W$ is an isometry, from \eqref{Ad 4} it follows that
		\begin{equation*}
			\begin{aligned}
				||p(T_1, \dots, T_{n-1}, T_n)||
				&= ||W^*p(M_{\Phi_1}, \dots, M_{\Phi_{n-1}}, M_z)W||\\
				&\le ||p(M_{\Phi_1}, \dots, M_{\Phi_{n-1}}, M_z)||\\
				&\le ||p||_{\infty, \overline{\mathbf{\Theta}}_n}.
			\end{aligned}
		\end{equation*}
		Therefore, $\mathbf{T}$ is a $\mathbf{\Theta}_n$-contraction. Also one can easily show that $\mathbf{T}^*$ is a $\mathbf{\Theta}_n$-contraction.
		
		Now it only remains to show that $(\widetilde{A}^{(i)}_0, \dots, \widetilde{A}^{(i)}_{p-1})$ is a solution of \eqref{Fundamental *}. For each $i$ we have that
		\begin{equation}\label{Ad 5}
			\begin{aligned}
				T^*_i - T_{n-i}T^{*p}_n
				&= W^*(M^*_{\Phi_i} - M_{\Phi_{n-i}}M^*p_z)W\\
				&= W^*\left(\sum_{l=0}^{p} M^{*l}_z \otimes \widetilde{A}^{(i)}_l - \left(\sum_{l=0}^{p} M^l_z \otimes \widetilde{A}^{(n-i)*}_l\right) \left(M^{*p}_z \otimes I_{\mathcal{D}_{T^*_n}}\right)\right)W.
			\end{aligned}
		\end{equation}
		We denote the projection onto the subspace spanned by $\{z^j : j \ge k\}$ and $\{1, z, \dots, z^j\}$ in $H^2(\mathcal{D}_{T^*_n})$ by $\mathbb{P}_{\overline{\langle z^k, z^{k+1}, \dots\rangle}}$ and $\mathbb{P}_{\overline{\langle 1, z, \dots, z^j \rangle}}$ respectively. Now, by the fact that $M^l_zM^{*l}_z = \mathbb{P}_{\overline{\langle z^k, z^{k+1}, \dots\rangle}}$ and simplifying \eqref{Ad 5}, we obtain
		\begin{equation}\label{Ad 6}
			\begin{aligned}
				T^*_i - T_{n-i}T^{*p}_n
				&= \sum_{j=0}^{p-1} W^*\mathbb{P}_{\overline{\langle 1, z, \dots, z^j \rangle}}\left(M^{*p-1-j}_z \otimes \widetilde{A}^{(i)}_{p-1-j}\right)W.
			\end{aligned}
		\end{equation}
		Each summand of \eqref{Ad 6} can be explicitly written as
		\begin{equation}\label{Ad 7}
			\begin{aligned}
				&W^*\mathbb{P}_{\overline{\langle 1, z, \dots, z^j \rangle}}\left(M^{*p-1-j}_z \otimes \widetilde{A}^{(i)}_{p-1-j}\right)Wh\\
				&= W^*\mathbb{P}_{\overline{\langle 1, z, \dots, z^j \rangle}}\left(M^{*p-1-j}_z \otimes \widetilde{A}^{(i)}_{p-1-j}\right)\left(\sum_{m \ge 0} z^m \otimes D_{T^*_n}T^{*m}h\right)\\
				&= W^*\mathbb{P}_{\overline{\langle 1, z, \dots, z^j \rangle}}\left(\sum_{m \ge p-1-j} z^{m-(p-1-j)} \otimes \widetilde{A}^{(i)}_{p-1-j}D_{T^*_n}T^{*m}h\right)\\
				&= W^*\left(\sum_{m=p-1-j}^{p-1} z^{m-(p-1-j)} \otimes \widetilde{A}^{(i)}_{p-1-j}D_{T^*_n}T^{*m}h\right)\\
				&= W^*\left(\sum_{m=0}^{j} z^m \otimes \widetilde{A}^{(i)}_{p-1-j}D_{T^*_n}T^{*m+p-1-j}h\right)\\
				&= \sum_{m=0}^{j} T^m_nD_{T^*_n} \widetilde{A}^{(i)}_{p-1-j}D_{T^*_n}T^{*m+p-1-j}h,
			\end{aligned}
		\end{equation}
		for all $h \in \mathcal{H}$. Then from \eqref{Ad 6} and \eqref{Ad 7} we have
		\begin{equation}\label{Ad 8}
			\begin{aligned}
				T^*_i - T_{n-i}T^{*p}_n
				&= \sum_{j-0}^{p-1} \sum_{m=0}^{j} T^m_n \widetilde{A}^{(i)}_{p-1-j}D_{T^*_n}T^{*m+p-1-j}.
			\end{aligned}
		\end{equation}
		By suitable rearrangement of the summands of \eqref{Ad 8} we obtain
		\begin{equation*}
			\begin{aligned}
				T^*_i - T_{n-i}T^{*p}_n
				&= \sum_{l=0}^{p-1} T^{p-1-l}_n D_{T^*_n}\widetilde{A}^{(i)}_lD_{T^*_n}T^{*p-1}_n + \sum_{l=0}^{p-2} T^{p-2-l}_n D_{T^*_n}\widetilde{A}^{(i)}_lD_{T^*_n}T^{*p-2}_n + \dots\\
				& \quad + \sum_{l=0}^{1} T^{1-l}_n D_{T^*_n}\widetilde{A}^{(i)}_lD_{T^*_n}T^*_n + D_{T^*_n}\widetilde{A}^{(i)}_0D_{T^*_n}, \quad 1 \le i \le n-1,
			\end{aligned}
		\end{equation*}
		which is as required. This finishes the proof.
	\end{proof}
	
	For a particular case, when $m = p = 1$, Theorem \ref{Admissible} produces only half of \cite[Theorem 3.4]{S. Pal}. That is why a very natural question could be asked: what about the fundamental operators of $\mathbf{T}$ in Theorem \ref{Admissible}? In \cite[Theorem 3.4]{S. Pal}, it is shown that the operator tuple $(A^{(1)}_0, \dots, A^{(n-1)}_0)$ is the fundamental operators of the $\Gamma_n$-contraction $\mathbf{T}$. In case of $m > 1, p > 1$, this is unresolved: what are the fundamental operators of $\mathbf{T}$? Therefore, the following problem could be interesting to work on.	
	
	\vspace{0.3cm}
	
	\noindent{\textbf{Open Question.}}
	Let $T_n$ be a pure contraction on a Hilbert space $\mathcal{H}$, and $\{A^{(i)}_l : 1 \le i \le n-1, 0 \le l \le p-1\}$ acting $\mathcal{D}_{T_n}$ and $\{\widetilde{A}^{(i)}_l : 1 \le i \le n-1, 0 \le l \le p-1\}$ satisfying \eqref{Compatibility Condition} and \eqref{Commutativity Condition}. Let
	\begin{equation*}
		\begin{aligned}
			\left(\gamma_1\left(\sum_{l=0}^{p} \widetilde{A}^{(1)*}_lz^l\right), \dots, \gamma_{n-1}\left(\sum_{l=0}^{p} \widetilde{A}^{(n-1)*}_lz^l\right)\right), \quad z \in \mathbb{T}
		\end{aligned}
	\end{equation*}
	is a $\Gamma_{n-1}$-contraction. Suppose $A^{(i)}_l$'s and $\widetilde{A}^{(i)}_l$'s satisfy \eqref{Admissible Condition}. Does there exists a $\mathbf{\Theta}_n$-contraction $\mathbf{T} = (T_1, \dots, T_n)$ such that $(A^{(i)}_0, \dots, A^{(i)}_{p-1})$ is a tuple of fundamental operators of $\mathbf{T}$ and $(\widetilde{A}^{(i)}_0, \dots, \widetilde{A}^{(i)}_{p-1})$ is a tuple of fundamental operators of $\mathbf{T}^*$?
	
	\vspace{0.3cm}
	
%
%
%
%

\end{document}